\documentclass[12pt,reqno]{amsart}

\usepackage[utf8]{inputenc}
\usepackage[english]{babel}
\usepackage{amsmath,amssymb,amsfonts,amsbsy,amsthm,amscd,latexsym,graphicx}
\usepackage{empheq,textcomp}

\theoremstyle{definition}
\newtheorem{theorem}{Theorem} [section]
\newtheorem{corollary}[theorem]{Corollary}
\newtheorem{lemma}[theorem]{Lemma}

\newtheorem{definition}[theorem]{Definition}

\newtheorem{remark}[theorem]{Remark}
\newtheorem{example}[theorem]{Example}
\numberwithin{equation}{section}

\newcommand{\Ac}{{\mathcal{A}}}

\newcommand{\Ec}{{\mathcal{E}}}
\newcommand{\Fc}{{\mathcal{F}}}
\newcommand{\Gc}{{\mathcal{G}}}
\newcommand{\Mcc}{\mathcal{M}}
\newcommand{\N}{\mathbb{N}}

\newcommand{\Yc}{\mathcal{Y}}

\newcommand{\Zc}{\mathcal{Z}}

\newcommand{\tFc}{{\widetilde{\mathcal{F}}}}
\newcommand{\tx}{{\widetilde{x}}}
\newcommand{\ty}{{\widetilde{y}}}
\newcommand{\tz}{{\widetilde{z}}}

\newcommand{\Eq}{\, = \,}

\newcommand{\Le}{\, \le \,}
\newcommand{\Lt}{\, < \,}
\newcommand{\Plus}{\, + \,}
\newcommand{\Minus}{\, - \,}
\newcommand{\Subseteq}{\, \subseteq \,}
\newcommand{\qeddef}{{\quad $\diamondsuit$}}
\newcommand{\qeddeff}{{\qquad \diamondsuit}}

\newcommand{\ip}[2]{\langle#1,#2\rangle}
\newcommand{\bigip}[2]{\bigl\langle #1, \, #2 \bigr\rangle}
\newcommand{\Bigip}[2]{\Bigl\langle #1, \, #2 \Bigr\rangle}
\newcommand{\biggip}[2]{\biggl\langle #1, \, #2 \biggr\rangle}
\newcommand{\norm}[1]{\|#1\|}
\newcommand{\Bignorm}[1]{\Bigl\|#1\Bigr\|}
\newcommand{\bigparen}[1]{\bigl(#1\bigr)}
\newcommand{\Bigparen}[1]{\Bigl(#1\Bigr)}
\newcommand{\biggparen}[1]{\biggl(#1\biggr)}
\newcommand{\set}[1]{\{#1\}}
\newcommand{\bigset}[1]{\bigl\{#1\bigr\}}
\newcommand{\Bigset}[1]{\Bigl\{#1\Bigr\}}
\newcommand{\clspan}{{\overline{\text{span}}}}
\newcommand{\takeaway}{\hskip 0.8 pt \backslash \hskip 0.8 pt}
\newcommand{\tinyspace}{\hskip 0.8 pt}

\begin{document}

\title{Convergence of Frame Series}

\author{Christopher Heil and Pu-Ting Yu}

\address{School of Mathematics, Georgia Institute of Technology,
Atlanta, GA 30332, USA}

\email{heil@math.gatech.edu}
\email{pyu73@gatech.edu}

\thanks{Acknowledgement:
This research was partially supported by a grant from the
Simons Foundation.}


\begin{abstract}
If $\set{x_n}_{n \in \N}$ is a frame for a Hilbert space $H,$
then there exists a canonical dual frame $\set{\tx_n}_{n \in \N}$
such that for every $x \in H$ we have $x = \sum \, \ip{x}{\tx_n} \, x_n,$
with unconditional convergence of this series.
However, if the frame is not a Riesz basis, then there exist
alternative duals $\set{y_n}_{n \in \N}$
and synthesis pseudo-duals $\set{z_n}_{n \in \N}$ such that
$x = \sum \, \ip{x}{y_n} \, x_n$
and $x = \sum \, \ip{x}{x_n} \, z_n$ for every $x.$
We characterize the frames for which the \emph{frame series}
\mbox{$x = \sum \, \ip{x}{y_n} \, x_n$}
converges unconditionally for every $x$ for every alternative dual,
and similarly for synthesis pseudo-dual.
In particular, we prove that if $\set{x_n}_{n \in \N}$
does not contain infinitely many zeros then the frame series
converge unconditionally for every alternative dual
(or synthesis pseudo-dual) if and only if
$\set{x_n}_{n \in \N}$ is a near-Riesz basis. 
We also prove that all alternative duals and synthesis pseudo-duals
have the same excess as their associated frame.
\end{abstract}

\maketitle

\section{Introduction}
Frames were first introduced by Duffin and Schaeffer in \cite{DS52}
in their study of non-harmonic Fourier series, and interest resurged
with the paper \cite{DGM86} by Daubechies, Grossmann, and Meyer,
which applied frames to wavelet and Gabor systems.
We refer to \cite{Chr16} and \cite{Hei11} for relatively recent
textbook recountings of the mathematics of frames.

We say that a sequence $\set{x_n}_{n \in \N}$
in a separable Hilbert space $H$ is a frame if there exist
positive constants $A \le B,$ called \emph{frame bounds}, such that 
$$A\, \norm{x}^2
\Le \sum_{n=1}^\infty |\ip{x}{x_n}|^2
\Le B\, \norm{x}^2,
\qquad\text{for all } x \in H.$$
A frame possesses basis-like properties, as
there exist sequences $\set{y_n}_{n \in \N}$
and $\set{z_n}_{n \in \N}$ such that
each $x$ in $H$ can be represented as
\begin{equation} \label{frameseries_eq}
x \Eq \sum_{n=1}^\infty \, \ip{x}{y_n} \, x_n 
\end{equation}
and
\begin{equation} \label{synthesis_frameseries_eq}
x \Eq \sum_{n=1}^\infty \, \ip{x}{x_n} \, z_n.
\end{equation}
Such a sequence $\set{y_n}_{n \in \N}$
(respectively $\set{z_n}_{n \in \N}$) is said to be an
\emph{alternative dual} (respectively \emph{synthesis pseudo-dual})
of $\set{x_n}_{n \in \N}.$
If in addition $\set{y_n}_{n \in \N}$ (respectively $\set{z_n}_{n \in \N}$)
is a frame then it is called an \emph{alternative dual frame}
(respectively \emph{synthesis pseudo-dual frame}).

In general, an alternative dual or synthesis pseudo-dual of a frame
need not be unique.
In some sense, frames exchange uniqueness of frame expansions
for flexibility in the choice of coefficients.
The \emph{canonical dual frame}, which we will define below,
is one alternative dual frame and it is also a synthesis pseudo-dual frame.
If $\set{y_n}_{n \in \N}$
is the canonical dual frame then the series in equation \eqref{frameseries_eq}
converges unconditionally for all $x \in H,$
and similarly for equation \eqref{synthesis_frameseries_eq}
if $\set{z_n}_{n \in \N}$ is the canonical dual frame.

The study of the convergence of frame series has a long history.
In \cite{Hei90}, Heil gave sufficient and necessary conditions
for a series $\sum c_n x_n$ to converge unconditionally,
where $\set{x_n}_{n \in \N}$ is a frame that is norm-bounded below.
In \cite{Hol94}, Holub established relations between types of frames
that he termed \emph{Besselian frames},
\emph{unconditional frames}, and \emph{near-Riesz bases}.
Convergence of Weyl--Heisenberg frame series was also studied in \cite{ZH06},
and unconditional constants were discussed in \cite{BCKL17}.
Stoeva and Balazs \cite{ST11, SB13, SB20}
studied the convergence of more general types of series
related to frames, called \textit{frame multipliers}.

We study the unconditional convergence of the series
$\sum \, \ip{x}{y_n} \, x_n$
and $\sum \, \ip{x}{x_n} \, z_n$
where $\set{y_n}_{n \in \N}$ is an alternative dual
and $\set{z_n}_{n \in \N}$ is a synthesis pseudo-dual
of $\set{x_n}_{n \in \N}.$
A natural question is whether this series must converge unconditionally
for every alternative dual and every $x,$
and similarly for synthesis pseudo-dual.
Barring the case where $\set{x_n}_{n \in \N}$ contains infinitely many zeros,
we prove that this holds true if and only if $\set{x_n}_{n \in \N}$
is a near-Riesz basis.
We also study the excess of alternative duals and synthesis pseudo-duals,
showing that the excess of these duals must be the same
as the excess of their associated frames.
 
Notation, terminology, and preliminary results will be presented
in Section~\ref{prelim_sec}.
In Section~\ref{converge_sec}, we address the convergence of
$\sum \, \ip{x}{y_n} \, x_n$ and $\sum \, \ip{x}{x_n} \, z_n,$
obtaining several equivalent characterizations of near-Riesz bases.
Finally, Section~\ref{excess_sec} studies the excess of alternative duals,
synthesis pseudo-duals, and their associated frames.

\smallskip
\section{Preliminaries} \label{prelim_sec}

Throughout this paper, $H$ denotes an infinite-dimensional
separable Hilbert space with inner product $\ip{\cdot}{\cdot}$
and corresponding norm $\norm{\cdot},$
and $\ell^2$ is the space of square-summable sequences
indexed by the natural numbers $\N.$

We list some basic results required for this paper, and refer to
\cite{Chr16} and \cite{Hei11} for more details and proofs.

If $\set{x_n}_{n \in \N}$ is a frame for $H,$ then its
\emph{frame operator} $Sx = \sum \, \ip{x}{x_n} \, x_n$
is a bounded linear invertible map of $H$ onto itself.
The \emph{canonical dual frame} is $\set{\tx_n}_{n \in \N},$
where $\tx_n = S^{-1}x_n.$
We have
$$x
\Eq \sum_{n=1}^\infty \, \ip{x}{\tx_n} \, x_n
\Eq \sum_{n=1}^\infty \, \ip{x}{x_n} \, \tx_n,
\qquad\text{for all } x \in H.$$
In particular, $\set{\tx_n}_{n \in \N}$ is an alternative dual frame
and also a synthesis pseudo-dual frame.

We say that $\set{x_n}_{n \in \N}$ is a \emph{Bessel sequence}
if it satisfies at least the upper inequality
in the definition of a frame.
That is, there exists some $B > 0,$
called an \emph{upper frame bound}, such that
$\sum |\ip{x}{x_n}|^2 \le B\, \norm{x}^2$
for every $x.$
In this case, 
the series $\sum c_n x_n$ converges unconditionally
whenever $(c_n)\in \ell^2.$

We recall the definition of a Riesz basis
and define unconditional frames and near-Riesz basis,
which were first introduced in \cite{Hol94}.
A variety of equivalent characterizations of Riesz bases
can be found in \cite[Thm.~8.32]{Hei11}.

\begin{definition} \label{nearriesz_def}
Let $\Fc = \set{x_n}_{n \in \N}$ be a frame for $H.$
\begin{enumerate}
\item[\textup{(a)}]
$\Fc$ is a \emph{Riesz basis} for $H$
if it is the image of an orthonormal basis for $H$
under a continuous linear bijection of $H$ onto itself.

\smallskip
\item[\textup{(b)}]
$\Fc$ is a \emph{near-Riesz basis} if there exists a
finite set $F \subseteq \N$ such that
$\set{x_n}_{n \notin F}$ is a Riesz basis.

\smallskip
\item[\textup{(c)}]
$\Fc$ is an \emph{unconditional frame}
if $\sum c_n x_n$ converges unconditionally whenever it converges. \qeddef
\end{enumerate}
\end{definition}

The notion of excess was introduced in \cite{BCHL03}
to address questions regarding overcompleteness of frames.
Letting $|\Gc|$ denote the cardinality of $\Gc$ if $\Gc$ is finite and
$|\Gc| = \infty$ otherwise, the
excess of a sequence $\Fc = \set{x_n}_{n \in \N}$ in $H$ is
$$e(\Fc)
\Eq \sup\bigset{|\Gc| : \Gc \subseteq \Fc \text{ and }
    \clspan(\Fc \backslash \Gc) = \clspan(\Fc)}.$$
A frame has finite excess if and only if it is a near-Riesz basis.
Holub \cite{Hol94} established that
unconditional frames and near-Riesz bases are equivalent
when the frame is norm-bounded below.
Casazza and Christensen \cite{CC96} generalized this to any frame
that does not contain infinitely many zeros.
We summarize those results as follows.

\begin{theorem} \label{nearriesz_thm}
If $\Fc = \set{x_n}_{n \in \N}$ is a sequence in $H,$
then the following two statements are equivalent.
\begin{enumerate}
\item[\textup{(a)}]
$\Fc$ is a near-Riesz basis.

\smallskip
\item[\textup{(b)}]
$\Fc$ is a frame that has finite excess.
\end{enumerate}

\smallskip\noindent
Moreover, if $\Fc$ does not contain infinitely many zeros then statements (a)
and (b) are also equivalent to the following statement.
\begin{enumerate}
\item[\textup{(c)}]
$\Fc$ is an unconditional frame. \qeddef
\end{enumerate}
\end{theorem}

The next result is Lemma 6.3.2 in \cite{Chr16}
(also see Theorem 8.18 in \cite{Hei11}).

\begin{lemma} \label{alternativedual_thm}
Let $\Fc = \set{x_n}_{n \in \N}$ be a frame for $H.$
\begin{enumerate}
\item[\textup{(a)}]
If an alternative dual $\set{y_n}_{n \in \N}$ is a Bessel sequence,
then $\set{y_n}_{n \in \N}$ is a frame for $H.$
In particular, an alternative dual frame
is also a synthesis pseudo-dual frame.

\smallskip
\item[\textup{(b)}]
If a synthesis pseudo-dual $\set{z_n}_{n \in \N}$ is a Bessel sequence,
then $\set{z_n}_{n \in \N}$ is a frame for $H.$
In particular, a synthesis pseudo-dual frame
is also an alternative dual frame. \qeddef
\end{enumerate}
\end{lemma}

The following is Theorem 2.2 in \cite{BT15}.

\begin{theorem} \label{alternativeexcess_thm}
If $\set{x_n}_{n \in \N}$ is a frame for $H$ and
$\set{y_n}_{n \in \N}$ is an alternative dual frame,
then $\set{x_n}_{n \in \N}$ and $\set{y_n}_{n \in \N}$
have the same excess. \qeddef
\end{theorem} 

It is known that, in at least some cases, convergence of the series
$\sum \, \ip{x}{y_n} \, x_n$ implies the convergence of
$\sum \, \ip{x}{x_n} \, y_n.$
For example, the following is a special case of \cite[Lem.~3.1]{SB13}.

\begin{lemma} \label{equiv_conv_alt_syn_dual}
Let $\set{x_n}_{n \in \N}$ and $\set{y_n}_{n \in \N}$
be two sequences in $H.$
Then $\sum \, \ip{x}{y_n} \, x_n$ converges unconditionally for every $x$
if and only if $\sum \, \ip{x}{x_n} \, y_n$ converges unconditionally
for every $x.$
In particular, if
$\sum \, \ip{x}{x_n} \, y_n$
and $\sum \, \ip{x}{y_n} \, x_n$
converge for every $x,$ then
$\sum \, \ip{x_0}{y_n} \, x_n$ converges conditionally for some $x_0 \in H$
if and only if $\sum \, \ip{y_0}{x_n} \, y_n$ converges conditionally
for some $y_0 \in H.$
\qeddef\end{lemma}

Finally, a sequence $\set{f_n}_{n \in \N}$ in $H$ is
\emph{minimal} if $f_m \notin \clspan\set{f_n}_{n \neq m}$
for every $m \in \N.$
This occurs if and only if there exists
a sequence $\set{g_n}_{n \in \N}$ that is \emph{biorthogonal} to
$\set{f_n}_{n \in \N}$; that is,
$\ip{f_m}{g_n} = \delta_{mn}$ for $m,$ $n \in \N$
(for one proof, see \cite[Lem.~5.4]{Hei11}).

\smallskip
\section{Convergence of Frame Series} \label{converge_sec}

If $\Fc = \set{x_n}_{n \in \N}$ is a frame then $\sum c_n x_n$
converges unconditionally whenever $(c_n)_{n \in \N}$ belongs to $\ell^2.$ 
Therefore
$\sum \, \ip{x}{y_n} \, x_n$ (respectively $\sum \, \ip{x}{x_n} \, y_n$)
converges unconditionally whenever $\set{y_n}_{n \in \N}$
is an alternative dual frame (respectively synthesis pseudo-dual frame).
Therefore, if every alternative dual (respectively synthesis pseudo-dual)
of $\Fc$ is a frame,
then unconditional convergence of $\sum \, \ip{x}{y_n} \, x_n$
(respectively $\sum \, \ip{x}{x_n} \, z_n$) is ensured.
However, we will construct a frame that possesses non-frame
alternative duals and non-frame synthesis pseudo-duals
for which $\sum \, \ip{u}{y_n} \, x_n$
and $\sum \, \ip{v}{x_n} \, z_n$
converge conditionally for some $u$ and $v.$
The frame $\Fc$ in this construction is norm-bounded below.

\begin{example} \label{dualnotframe_example}
Let $\Ec = \set{e_n}_{n \in \N}$ be an orthonormal basis for $H,$
and consider the frame
$$\Fc
\Eq \set{x_n}_{n \in \N}
\Eq \set{e_1,e_1,e_1,e_2,e_2,e_3,e_3,e_4,e_4,\dots}.$$
Define $y_1=e_1,$ $y_2=e_1,$ and $y_3=-e_1,$ and for $n \geq 2$ set 
$$y_{2n} \Eq \frac{1}{\sqrt{n}} \, e_1 \Plus \frac{1}{2} \, e_n
\qquad\text{and}\qquad
y_{2n+1} \Eq -\frac{1}{\sqrt{n}} \, e_1 \Plus \frac{1}{2} \, e_n.$$
If $x \in H,$ then for each integer $K \geq 4$ we have that
\begin{equation*}
\sum_{n=1}^K \, \ip{x}{y_n} \, x_n 
\Eq \begin{cases}
    \displaystyle \sum^{(K-1)/2}_{n=1} \ip{x}{e_n} \, e_n,
    & \text{if $K$ is odd,} \\[5 \jot]
    \displaystyle \sum_{n=1}^{(K/2)-1} \ip{x}{e_n} \, e_n \Plus
    \sqrt{\frac{2}{K}} \, \ip{x}{e_1} \, e_{K/2} \Plus
    \frac{1}{2} \, \ip{x}{e_{K/2}} \, e_{K/2},
    & \text{if $K$ is even},
\end{cases}
\end{equation*}
and 
\begin{equation*}
\sum_{n=1}^K \, \ip{x}{x_n} \, y_n 
\Eq \begin{cases}
    \displaystyle \sum^{(K-1)/2}_{n=1} \ip{x}{e_n} \, e_n,
    & \text{if $K$ is odd,} \\[5 \jot]
    \displaystyle \sum_{n=1}^{(K/2)-1} \ip{x}{e_n} \, e_n \Plus
    \sqrt{\frac{2}{K}} \, \ip{x}{e_{K/2}} \, e_{1} \Plus
    \frac{1}{2} \, \ip{x}{e_{K/2}} \, e_{K/2},
    & \text{if $K$ is even.}
\end{cases}
\end{equation*}
Since both $\sqrt{2/K}$ and $\ip{x}{e_{K/2}}$ tend to zero as $K$ increases,
it follows that
\begin{equation} \label{dualexample_eq}
\sum_{n=1}^\infty \, \ip{x}{y_n} \, x_n
\Eq x
\Eq \sum_{n=1}^\infty \, \ip{x}{x_n} \, y_n.
\end{equation}
Therefore $\Yc = \set{y_n}_{n \in \N}$ is both an alternative dual
and a synthesis pseudo-dual of $\Fc.$
However, for $x = e_1$ the first representation
$e_1 = \sum_{n=1}^\infty \, \ip{e_1}{y_n} \, x_n$
from equation \eqref{dualexample_eq} is
$$e_1
\Eq e_1 \Plus e_1 \Minus e_1 \Plus \frac{1}{\sqrt{2}}\, e_2
    \Minus \frac{1}{\sqrt{2}} \, e_2 \Plus \cdots \Plus \frac{1}{\sqrt{n}} \, e_n
    \Minus \frac{1}{\sqrt{n}} \, e_n \Plus \cdots,$$
which converges conditionally.
Consequently $\Yc$ cannot be a frame. 
Moreover, Lemma \ref{equiv_conv_alt_syn_dual} implies that
there exists some $v$ such that the second representation
$v = \sum_{n=1}^\infty \, \ip{v}{x_n} \, y_n$
in equation \eqref{dualexample_eq} must converges conditionally.
\qeddef\end{example}

Several different examples illustrating convergence of frame series
are given in \cite{ST11}.

The frame $\Fc$ in Example \ref{dualnotframe_example}
has infinite excess.
We prove next that if a frame $\set{x_n}_{n \in \N}$ has finite excess,
then for every $x$ the frame series $\sum \, \ip{x}{y_n} \, x_n$
and $\sum \, \ip{x}{x_n} \, z_n$ converge unconditionally
for every alternative dual $\set{y_n}_{n \in \N}$
and synthesis pseudo-dual $\set{z_n}_{n \in \N}.$
We remark that a different proof
can be given by using \cite[Prop.~5.11]{ST11}.

\begin{theorem} \label{finiteexcess_thm}
If $\Fc = \set{x_n}_{n \in \N}$ is a near-Riesz basis for $H,$
then every alternative dual and synthesis pseudo-dual of $\Fc$ 
are near-Riesz bases and hence are frames.
\end{theorem}
\begin{proof}
By Lemma \ref{alternativedual_thm} and Theorem \ref{alternativeexcess_thm},
it suffices to show that every alternative dual
$\Yc = \set{y_n}_{n \in \N}$ and every
synthesis pseudo-dual $\Zc = \set{z_n}_{n \in \N}$ is a Bessel sequence.

Let $A = \set{n_j}_{j=1}^N$ be a subset of $\N$
such that $\set{x_n}_{n \notin A}$ is a Riesz basis for $H.$
Then $\set{x_n}_{n \notin A}$ has a biorthogonal sequence
$\set{\tx_n}_{n \notin A},$ and consequently
\begin{equation} \label{notinA_eq}
x \Eq \sum_{n \notin A} \, \ip{x}{\tx_n} \, x_n
\end{equation}
is the unique representation of $x$ in terms of $\set{x_n}_{n \notin A}.$

Fix $x \in H.$
Since $\set{y_n}_{n \in \N}$ is an alternative dual of
$\set{x_n}_{n \in \N},$
\begin{equation} \label{altdual_eq}
x
\Eq \sum_{n=1}^\infty \, \ip{x}{y_n} \, x_n
\Eq \sum_{n \notin A} \, \ip{x}{y_n} \, x_n \Plus
    \sum_{j=1 }^N \, \ip{x}{y_{n_j}} \, x_{n_j}.
\end{equation}
By equation \eqref{notinA_eq}, for each $j = 1, \dots, N$ we have
$x_{n_j}
= \sum_{n \notin A} \, \ip{x_{n_j}}{\tx_n} \, x_n.$
Substituting this into equation \eqref{altdual_eq} gives
\begin{equation} \label{newnotinA_eq}
x
\Eq \sum_{n \notin A} \, \biggparen{\ip{x}{y_n} \Plus
    \sum_{j=1}^N \, \ip{x}{y_{n_j}} \, \ip{x_{n_j}}{\tx_n}} \, x_n.
\end{equation}
Since equation \eqref{notinA_eq} is the unique representation of $x$
in terms of $\set{x_n}_{n \notin A},$ we must therefore have
$$\ip{x}{y_n} \Plus \sum_{j=1}^N \, \ip{x}{y_{n_j}} \, \ip{x_{n_j}}{\tx_n}
\Eq \ip{x}{\tx_n}
\qquad\text{for } n \notin A.$$
Consequently, as sequences in $\ell^2(\N \takeaway A),$
$$\bigset{\ip{x}{y_n}}_{n \notin A}
\Eq \bigset{\ip{x}{\tx_n}}_{n \notin A}
    \Minus \sum_{j=1}^N \bigset{\ip{x}{y_{n_j}} \,
    \ip{x_{n_j}}{\tx_n}}_{n \notin A}.$$
Let $L$ denote an upper frame bound for
$\set{\tx_n}_{n \notin A}.$
Then by applying the Triangle and
Cauchy--Bunyakowski--Schwarz Inequalities, we compute that
\begin{align*}
\Bignorm{\bigset{\ip{x}{y_n}}_{n \notin A}}_{\ell^2(\N \takeaway A)}
& \Le \Bignorm{\bigset{\ip{x}{\tx_n}}_{n \notin A}}_{\ell^2(\N \takeaway A)}
      \Plus \sum_{j=1}^N \, \Bignorm{\bigset{\ip{x}{y_{n_j}} \,
      \ip{x_{n_j}}{\tx_n}}_{n \notin A}}_{\ell^2(\N \takeaway A)}
      \\
& \Le L^{1/2} \, \norm{x} \Plus \sum_{j=1}^N \, \norm{x} \, \norm{y_{n_j}} \,
      \Bignorm{\bigset{\ip{x_{n_j}}{\tx_n}}_{n \notin A}}
      _{\ell^2(\N \takeaway A)}
      \allowdisplaybreaks \\[0.5 \jot]
& \Le L^{1/2} \, \norm{x} \Plus \sum_{j=1}^N \, \norm{x} \, \norm{y_{n_j}} \,
      L^{1/2} \, \norm{x_{n_j}}
      \\[0.5 \jot]
& \Eq K \, \norm{x},
\end{align*}
where $K$ is a constant independent of $x.$
Therefore
\begin{align*}
\sum_{n=1}^\infty |\ip{x}{y_n}|^2
& \Eq \sum_{n \notin A} |\ip{x}{y_n}|^2 \Plus \sum_{n \in A} |\ip{x}{y_n}|^2
\Le K^2 \, \norm{x}^2 \Plus \sum_{n \in A} \, \norm{x}^2 \, \norm{y_n}^2
\Eq C \, \norm{x}^2,
\end{align*}
where $C$ is a constant independent of $x.$
Thus $\set{y_n}_{n \in \N}$ is a Bessel sequence in $H.$

On the other hand, if $k \notin A$ then
\begin{equation} \label{syn_dual_eq}
\tx_k
\Eq \sum_{n=1}^\infty \, \ip{\tx_k}{x_n} \, z_n
\Eq z_k\Plus
\sum_{j=1}^N \, \ip{\tx_k}{x_{n_j}} \, z_{n_j}.
\end{equation}
Consequently, if $x\in H$ then
$$\bigset{\ip{x}{z_k}}_{k \notin A}
\Eq \bigset{\ip{x}{\tx_k}}_{k \notin A}
    \Minus \sum_{j=1}^N \bigset{\ip{x}{z_{n_j}} \,
    \ip{x_{n_j}}{\tx_k}}_{k \notin A}.$$
Arguing as before, it follows that $\set{z_n}_{n \in \N}$
is a Bessel sequence.
\end{proof}

\begin{corollary} \label{nearriesz_cor}
If $\set{x_n}_{n \in \N}$ is a near-Riesz basis for $H,$ then
$\sum \, \ip{x}{y_n} \, x_n$ (respectively $\sum \, \ip{x}{x_n} \, z_n$)
converges unconditionally for every $x \in H$ and every alternative dual
$\set{y_n}_{n \in \N}$
(respectively synthesis pseudo-dual $\set{z_n}_{n \in \N}$).
\qeddef\end{corollary}

We will prove a converse to Theorem \ref{finiteexcess_thm} and
Corollary \ref{nearriesz_cor} below
(see Theorem \ref{nearrieszcharacterization_thm}).
To this end, we first prove that if $c_n$ are scalars such that the series
$x = \sum c_n x_n$ converges to a nonzero element in $H,$ then there exists
an alternative dual $\set{y_n}_{n \in \N}$ such that $c_n = \ip{x}{y_n}$
for every~$n.$

\begin{theorem} \label{realizable_thm}
Let $\Fc = \set{x_n}_{n \in \N}$ be a frame for $H$ and let
$(c_n)_{n \in \N}$ be a sequence of scalars.
If $\sum c_n x_n$ converges to a nonzero vector $x_0 \in H,$
then there exists a sequence $\set{y_n}_{n \in \N}$ that is
both an alternative dual and a synthesis pseudo-dual of $\Fc$
such that $c_n = \ip{x_0}{y_n}$ for every~$n.$
\end{theorem}
\begin{proof}
Assume that $\sum c_n x_n = x_0 \neq 0$
and let $\set{\tx_n}_{n \in \N}$ be the canonical dual frame of $\Fc.$
Set $M = \clspan\set{x_0},$ and let $P_{M^\perp}$ denote the
orthogonal projection of $H$ onto $M^\perp.$
Define
$$y_n
\Eq \frac{\overline{c_n}}{\norm{x_0}^2} \, x_0 \Plus P_{M^\perp} \tx_n,
\qquad\text{for } n \in \N.$$
We will show that $\set{y_n}_{n \in \N}$ has the required properties.

Fix $x \in H,$ and write $x$ as $x = x_M + x_{M^\perp}$
where $x_M \in M$ and $x_{M^\perp} \in M^\perp.$
Since $M$ is one-dimensional, there exists a scalar $\alpha$
such that $x_M = \alpha x_0.$
Therefore,
\begin{align*}
\sum_{n=1}^\infty \ip{x}{y_n} \, x_n
& \Eq \sum_{n=1}^\infty \, \Bigip{x_M \Plus x_{M^\perp}}
      {\frac{\overline{c_n}}{\norm{x_0}^2} \, x_0
      \Plus P_{M^\perp}\tx_n} \, x_n
      \\[0.5 \jot]
& \Eq \sum_{n=1}^\infty \,
      \Bigip{\alpha x_0}{\frac{\overline{c_n}}{\norm{x_0}^2} \, x_0} \, x_n
      \Plus \sum_{n=1}^\infty \, \bigip{x_{M^\perp}}{P_{M^\perp}\tx_n} \, x_n
      \qquad\text{\small (cross terms vanish)}
      \allowdisplaybreaks \\[0.5 \jot]
& \Eq \sum_{n=1}^\infty \alpha \tinyspace c_n x_n \Plus \sum_{n=1}^\infty \,
      \bigip{x_{M^\perp}}{P_{M^\perp}\tx_n} \, x_n
      \\
& \Eq \alpha x_0 \Plus \sum_{n=1}^\infty \, \ip{x_{M^\perp}}{\tx_n} \, x_n
\Eq x_M \Plus x_{M^\perp}
\Eq x.
\end{align*} 
This shows that $\set{y_n}_{n \in \N}$ is an alternative dual of $\Fc.$
Further, for every $n$ we have
$$\ip{x_0}{y_n}
\Eq \Bigip{x_0}{\frac{\overline{c_n}}
    {\norm{x_0}^2} \, x_0 \Plus P_{M^\perp}\tx_n}
\Eq c_n.$$

In order to show that $\set{y_n}_{n \in \N}$
is a synthesis pseudo-dual of $\set{x_n}_{n \in \N},$
fix any $x \in H.$
Then, since $\set{\tx_n}_{n \in \N}$ is a frame,
Lemma \ref{alternativedual_thm} implies that
$$\sum_{n=1}^\infty \, \ip{x}{x_n} \, P_{M^\perp}\tx_n
\Eq P_{M^\perp}\biggparen{\sum_{n=1}^\infty \, \ip{x}{x_n} \, \tx_n}
\Eq P_{M^\perp}\biggparen{\sum_{n=1}^\infty \, \ip{x}{\tx_n} \, x_n}
\Eq x_{M^\perp}.$$
On the other hand,
$$\sum_{n=1}^\infty \, \ip{x}{x_n}\, \frac{\overline{c_n}}{\norm{x_0}^2} \, x_0
\Eq \biggip{x}{\sum_{n=1}^\infty \,\frac{c_n}{\norm{x_0} \, }x_n} \,
    \frac{x_0}{\norm{x_0}}
\Eq \biggip{x}{\frac{x_0}{\norm{x_0}}}\, \frac{x_0}{\norm{x_0}}
\Eq x_M,$$
the last equality following from the fact that $x_0/\norm{x_0}$
is a unit vector.
Consequently,
$$\sum_{n=1}^\infty \, \ip{x}{x_n} \,y_n
\Eq \sum_{n=1}^\infty \, \ip{x}{x_n} \,
    \biggparen{\frac{\overline{c_n}}{\|x_0\|^2}\,x_0+P_{M^\perp}\tx_n}
\Eq x_M \Plus x_{M^\perp}
\Eq x. \qquad\qedhere$$
\end{proof}

We obtain some corollaries of Theorem \ref{realizable_thm}.

\begin{definition} \label{realizable_def}
Let $\Fc = \set{x_n}_{n \in \N}$ be a frame for $H.$
We say that a sequence of scalars $(c_n)_{n \in \N}$
is \emph{realizable} with respect to $\Fc$
if there exists an $x \in H$ and an alternative dual
$\set{y_n}_{n \in \N}$ such that $c_n = \ip{x}{y_n}$ for all $n.$
\qeddef\end{definition}

\begin{corollary} \label{realizable_cor}
Let $\Fc = \set{x_n}_{n \in \N}$ be a frame for $H.$
Then a sequence $(c_n)_{n \in \N}$ is realizable with respect to $\Fc$
if and only if either $\sum c_n x_n$ converges to a nonzero $x \in H$
or $c_n = 0$ for every $n.$
\end{corollary}
\begin{proof}
$(\Rightarrow)$
Assume that $(c_n)_{n \in \N}$ is realizable.
Then there is some $x \in H$ and alternative dual $\set{y_n}_{n \in \N}$
such that $c_n = \ip{x}{y_n}$ for every $n.$
Then $x = \sum \, \ip{x}{y_n} \, x_n = \sum c_n x_n.$
If $x = 0,$ then $c_n = \ip{x}{y_n} = 0$ for every $n.$

\smallskip
$(\Leftarrow)$
If $c_n = 0$ for every $n$ then take $x = 0$ and let
$\set{y_n}_{n \in \N}$ be any alternative dual.
The case $\sum c_n x_n = x_0 \neq 0$ follows from
Theorem \ref{realizable_thm}.
\end{proof}

\begin{corollary} \label{conditional_cor}
Let $\Fc = \set{x_n}_{n \in \N}$ be a frame for $H.$
If there exist scalars $c_n$ such that $\sum c_n x_n$ converges conditionally,
then there exists an alternative dual $\set{y_n}_{n \in \N}$
that is also a synthesis pseudo-dual
and some vectors $x,$ $y \in H$ such that the series
$\sum \, \ip{x}{y_n}\, x_n$ and $\sum \, \ip{y}{x_n}\, y_n$
converge conditionally.
\end{corollary}
\begin{proof}
If $\sum c_n x_n \neq 0,$ then this follows from
Theorem \ref{realizable_thm} and Lemma \ref{equiv_conv_alt_syn_dual}.
Therefore, suppose that $\sum c_n x_n$ converges to zero conditionally.
Then $c_{n_0} x_{n_0} \neq 0$ for some $n_0 \in \N.$
Let $d_n = c_n$ for $n \neq n_0$ and set $d_{n_0} = 2c_{n_0}.$
Then $\sum d_n x_n$ converges conditionally to $c_{n_0} x_{n_0},$
which is not zero.
Hence we can apply Theorem \ref{realizable_thm}
and Lemma \ref{equiv_conv_alt_syn_dual} again.
\end{proof}

\begin{corollary} \label{nearrieszcharacterization_cor}
Assume that $\Fc = \set{x_n}_{n \in \N}$ is a frame for $H$
that contains at most finitely many zeros.
Then the following statements are equivalent.
\begin{enumerate}
\item[\textup{(a)}]
$\Fc$ is a near-Riesz basis.

\smallskip
\item[\textup{(b)}]
$\sum \, \ip{x}{y_n} \, x_n$ converges unconditionally
for all alternative duals $\set{y_n}_{n \in \N}$ and every $x \in H.$

\smallskip
\item[\textup{(c)}]
$\sum \, \ip{x}{x_n} \, z_n$ converges unconditionally
for all synthesis pseudo-duals $\set{z_n}_{n \in \N}$ and every $x \in H.$
\end{enumerate}
\end{corollary}
\begin{proof}
The implications (a) $\Rightarrow$ (b) and (a) $\Rightarrow$ (c)
follow from Theorem \ref{finiteexcess_thm}.

For the converse directions, suppose that $\Fc$ is not a near-Riesz basis.
Then, by Theorem \ref{nearriesz_thm}, $\Fc$ is not an unconditional frame.
Therefore there are some scalars $c_n$ such that $\sum c_n x_n$
converges conditionally.
The result therefore follows by Corollary \ref{conditional_cor}.
\end{proof}

\begin{remark}
The key point in the proof of the implication
(b) $\Rightarrow$ (a) in Corollary \ref{nearrieszcharacterization_cor}
is that we can realize $c_n$ as $\ip{x}{y_n}$
for some alternative dual $\set{y_n}_{n \in \N}$ and $x \in H$
if $\sum \, c_n x_n$ converges to some nonzero element.
We can generalize this as follows.
The proof of Theorem \ref{realizable_thm} for the case of
alternative duals does not utilize all of the properties of frame.
In fact, we do not need the existence of a canonical dual;
instead, we only require the existence of at least one sequence
$\set{\ty_n}_{n \in \N}$
such that $x = \sum \, \ip{x}{\ty_n}\, x_n$ for every $x.$
Hence, whenever sequences $\set{x_n}_{n \in \N}$ and
$\set{\ty_n}_{n \in \N}$ have this property, it will be the case that
a sequence of scalars $(c_n)_{n \in \N}$ can be expressed as
$\set{\ip{x}{y_n}}_{n \in \N}$ for some $x$ and alternative dual
$\set{y_n}_{n \in \N}$
if $\sum \, c_n x_n$ converges to a nonzero element.
Here, ``alternative dual'' simply means a sequence such that
$x = \sum \, \ip{x}{y_n} \, x_n$ for every $x \in H$. 

Based on this observation, we can give a more general characterization
of sequences for which $x = \sum \, \ip{x}{y_n} \, x_n$
converges unconditionally for every $x$ and
alternative dual $\set{y_n}_{n \in \N}$.
In \cite{CC96}, Casazza and Christensen showed that if a sequence
$\set{x_n}_{n \in \N}$ containing at most finitely many zeros is such that
$\sum c_n x_n$ converges unconditionally whenever it converges, then
$\set{x_n}_{n \in \N}$ is an unconditional basis plus at most
finitely many elements.
Therefore, if $\set{x_n}_{n \in \N}$ is a sequence that
contains at most finitely many zeros and it is the case
that whenever $\set{y_n}_{n \in \N}$ is a sequence
such that $x = \sum \, \ip{x}{y_n} \, x_n$ for every $x$
then this series converges unconditionally for every $x,$
then a proof similar to that of Corollary \ref{nearrieszcharacterization_cor}
shows that $\set{x_n}_{n \in \N}$ must be
an unconditional basis plus at most finitely many elements.
\qeddef\end{remark}

Now we prove that the converse of Theorem \ref{finiteexcess_thm} holds.

\begin{theorem} \label{nearrieszcharacterization_thm}
Let $\Fc = \set{x_n}_{n \in \N}$ be a frame for $H.$
Then the following statements are equivalent.
\begin{enumerate}
\item[\textup{(a)}]
$\Fc$ is a near-Riesz basis.

\smallskip
\item[\textup{(b)}]
Every alternative dual of $\Fc$ is a frame.

\smallskip
\item[\textup{(c)}]
Every synthesis pseudo-dual of $\Fc$ is a frame.
\end{enumerate}
\end{theorem}
\begin{proof}
The implication (a) $\Rightarrow$ (b) was proved
in Theorem \ref{finiteexcess_thm}.

For the implication (b) $\Rightarrow$ (a),
assume that $\Fc$ is not a near-Riesz basis.
We must show that there is an alternative dual of $\Fc$ that is not a frame.

First consider the case where $\Fc$ contains infinitely many zeros.
Since it is a frame it must also contain infinitely many nonzero elements,
so by reindexing we may assume that $x_{2n} = 0$ for every $n.$
Let $\set{y_n}_{n \in \N}$ be any alternative dual.
Choose any $x_0 \neq 0$ in $H,$ and define a new sequence
$\set{w_n}_{n \in \N}$ by
$w_{2n-1} = y_{2n-1}$ and $w_{2n} = x_0.$
Then $\set{w_n}_{n \in \N}$ is an alternative dual of $\Fc,$
but it is not a frame.

On the other hand, if $\Fc$ does not contain infinitely many zeros,
then Corollary \ref{nearrieszcharacterization_cor} implies that
there exists an alternative dual $\set{y_n}_{n \in \N}$ and some $x \in H$
such that $\sum \, \ip{x}{y_n} \, x_n$ converges conditionally.
But then $\set{y_n}_{n \in \N}$ is not a frame. 

A similar argument shows the equivalence of (a) and (c).
\end{proof}

\begin{remark}
We have seen that excess of a frame is related to the
unconditional convergence of frame series.
Consequently, if we use unconditional convergence of frame series
as a criterion to distinguish ``good'' and ``bad'' frames,
then only near-Riesz bases can be ``good'' frames.
A natural follow-up question is:
For which alternative dual $\set{y_n}_{n \in \N}$
(respectively synthesis pseudo-dual)
can the corresponding frame expansion
$\sum \, \ip{x}{y_n} \, x_n$
(respectively $\sum \, \ip{x}{x_n} \, y_n$)
converge unconditionally for every $x \in H$?
If the frame $\set{x_n}_{n \in \N}$ is norm-bounded below then
we know that an alternative dual
(respectively synthesis pseudo-dual)
can be ``good'' in this sense if it is a frame.
For two arbitrary sequences $\set{x_n}_{n \in \N}$ and $\set{y_n}_{n \in \N}$
in $H,$ it was conjectured in \cite{SB13} that the series
$\sum \, \ip{x}{y_n} \, x_n$
converges unconditionally for every $x \in H$ if and only if
there exist sequences of scalars
$(c_n)_{n \in \N}$ and $(d_n)_{n \in \N}$ such that
$c_n \overline{d_n} =1 $ for every $n$ and
$\set{c_n x_n}_{n \in \N}$ and $\set{d_n y_n}_{n \in \N}$
are Bessel sequences.
\qeddef\end{remark}

A frame is a Riesz basis if and only if the range of the analysis operator
$Cx = \set{\ip{x}{x_n}}_{n \in \N}$
is $\ell^2$ (see \cite[Thm.~8.32]{Hei11}).
We know that a frame that is not a Riesz basis possesses
more than one alternative dual, in fact it has infinitely many
(see \cite[Lem.~6.3.1]{Chr16}).
Further, in this case the range of the analysis operator
is a proper closed subspace of $\ell^2.$
However, could it be that the union of the ranges of the analysis operators
of all alternative duals is $\ell^2$?
We address this question next.

\begin{definition} \label{momentspace_def}
(a) The \emph{moment space} associated with a sequence
$\Fc = \set{f_n}_{n \in \N}$ in $H$ is
$$m(\Fc)
\Eq \Bigset{\set{\ip{x}{f_n}}_{n \in \N} : x \in H}.$$

(b) The \emph{extended moment space} of a frame
$\Fc = \set{f_n}_{n \in \N}$ for $H$ is the union of all
moment spaces over all alternative duals:
$$\Mcc(\Fc)
\Eq \bigcup \, \bigset{m(\Yc) : \Yc \text{ is an alternative dual of } \Fc}.
    \qeddeff$$
\end{definition}

For details on moment spaces, we refer to \cite{You01}.
The following result is Theorem~7 in Chapter~4 of \cite{You01}.

\begin{lemma} \label{moment_lemma}
If $\Fc = \set{f_n}_{n \in \N}$ and $\Gc = \set{g_n}_{n \in \N}$
are two complete sequences in $H,$ then $m(\Fc) = m(\Gc)$
if and only if there exists an bounded linear invertible operator
$T \colon H \to H$ such that $T(f_n) = g_n$ for every $n \in \N.$
\qeddef\end{lemma}

We will need the following lemma.

\begin{lemma} \label{momentperp_lemma}
Let $\Fc = \set{x_n}_{n \in \N}$ be a frame for $H,$ and let
$\tFc = \set{\tx_n}_{n \in \N}$ be its canonical dual frame.
Then for any alternative dual or synthesis pseudo-dual
$\set{y_n}_{n \in \N}$ of $\Fc,$
$$\set{\ip{y}{y_n}}_{n \in \N} \in m(\tFc\,)^\perp \ \iff \ y = 0.$$
\end{lemma}
\begin{proof}
Assume that $\set{\ip{y}{y_n}}_{n \in \N}$ belongs to $m(\tFc\,)^\perp.$
The frame operator $S$ for $\set{x_n}_{n \in \N}$ is a
bounded linear invertible map of $H$ onto itself and the
canonical dual frame is given by $\tx_n = S^{-1}x_n.$
Therefore Lemma \ref{moment_lemma} implies that $m(\tFc\,) = m(\Fc).$
Choose any $x \in H.$
Then $\set{\ip{x}{x_n}}_{n \in \N}$ belongs to $m(\Fc),$
so if $\set{y_n}_{n \in \N}$ is an alternative dual then
\begin{align*}
\ip{x}{y}
& \Eq \biggip{x}{\sum_{n=1}^\infty \, \ip{y}{y_n} \, x_n}
\Eq \sum_{n=1}^\infty \, \ip{x}{x_n} \, \overline{\ip{y}{y_n}}
\Eq 0.
\end{align*} 
Similarly, if $\set{y_n}_{n \in \N}$ is a synthesis pseudo-dual, then
\begin{align*}
\ip{x}{y}
& \Eq \biggip{\sum_{n=1}^\infty \, \ip{x}{x_n} \, y_n}{y}
\Eq \sum_{n=1}^\infty \, \ip{x}{x_n} \, \overline{\ip{y}{y_n}}
\Eq 0.
\end{align*} 
In either case we see that $\ip{x}{y} = 0$ for every $x,$ so $y = 0.$
\end{proof}

We now characterize the relation between $\Mcc(\Fc)$ and $\ell^2.$

\begin{theorem} \label{extendedmoment_thm}
If $\Fc = \set{x_n}_{n \in \N}$ is a frame for $H,$
then $\Mcc(\Fc) \subseteq \ell^2$ if and only if $\Fc$ is a near-Riesz basis.
Consequently, $\Mcc(\Fc) = \ell^2$ if and only if $\Fc$ is a Riesz basis.
\end{theorem}
\begin{proof}
By Theorem \ref{finiteexcess_thm}, if $\Fc$ is a near-Riesz basis
then $\Mcc(\Fc) \subseteq \ell^2$

Conversely, if $\Mcc(\Fc) \subseteq \ell^2$ then every alternative dual
is a Bessel sequence and consequently a frame by
Theorem \ref{alternativedual_thm}.
Therefore $\Fc$ is a near-Riesz basis by
Theorem \ref{nearrieszcharacterization_thm}.

Now assume that $\Fc$ is a Riesz basis.
In this case the canonical dual frame is also a Riesz basis,
so $\Mcc(\Fc) = \ell^2.$

Conversely, suppose that $\Mcc(\Fc) = \ell^2.$
Let $\tFc = \set{\tx_n}_{n \in \N}$ be the canonical dual frame of $\Fc,$
but suppose that $m(\tFc\,) \ne \ell^2.$
Then $m(\tFc\,)^\perp$ contains a nonzero sequence,
so since $\Mcc(\Fc) = \ell^2$ there is some nonzero $y \in H$ such that
$\set{\ip{y}{y_n}}_{n \in \N} \in m(\tFc\,)^\perp.$
But then $y=0$ by Lemma \ref{momentperp_lemma},
which is a contradiction.
Therefore $m(\tFc\,) = \ell^2,$ which implies that $\tFc,$ and hence $\Fc,$
is a Riesz basis.
\end{proof}

We can also define the extended moment space for synthesis pseudo-duals.
Since the canonical dual frame is a synthesis pseudo-dual frame,
we can use a similar proof to obtain a synthesis pseudo-dual version
of Theorem \ref{extendedmoment_thm}.

\smallskip
\section{The Excess of Alternative Duals and Synthesis Pseudo-Duals}
\label{excess_sec}

Now we consider the relation between the excess of a frame and the
excess of its alternative duals and synthesis pseudo-duals.
We know from Theorem \ref{alternativeexcess_thm}
that every alternative dual that is a frame has the same excess
as its associated frame.
Moreover, we proved in Theorem \ref{finiteexcess_thm} that
every alternative dual (respectively synthesis pseudo-dual)
of a frame with finite excess must be a frame.
In particular, every alternative dual frame is also a
synthesis pseudo-dual frame by Lemma \ref{alternativedual_thm}.
Consequently, if $\Fc$ is a frame with finite excess
then $e(\Fc) = e(\Yc)$ for every alternative dual and
synthesis pseudo-dual $\Yc.$
We will prove in this section that this same relation
holds for frames with infinite excess.

First we need the following lemma.

\begin{lemma} \label{minimal_lemma}
If $\Fc = \set{f_n}_{n \in \N}$ is a minimal sequence in $H$ and
$\Gc = \set{g_n}_{n=1}^N$ is a finite sequence, then
$e(\Fc \cup \Gc) \Lt \infty.$
\end{lemma}
\begin{proof}
Observe that if $J$ is a finite subset of $\N,$
then $\text{codim}\bigparen{\clspan\set{f_n}_{n \notin J}} \geq |J|$
since $\set{f_n}_{n \in \N}$ is minimal.

By replacing $H$ with the closed span of $\Fc \cup \Gc,$
we may assume that $\Fc \cup \Gc$ is complete.
Suppose that $e(\Fc \cup \Gc) = \infty.$
Then we can find $3N$ elements of $\Fc \cup \Gc$ such that
the removal of these $3N$ elements still leaves a complete sequence.
Precisely, there exist sets
$J_1 \subseteq \N$ and $J_2 \subseteq \set{1,\dots,N}$
such that $|J_1| + |J_2| = 3N$ and
$$\Ac
\Eq \set{f_n}_{n \in \N \takeaway J_1}
    \cup \set{g_n}_{n \in \set{1,\dots,N} \takeaway J_2}
    \ \text{ is complete}.$$
Since $|J_2| \leq N,$ we must have $|J_1|\geq 2N.$
Therefore
$\text{codim}\bigparen{\clspan\set{f_n}_{n \in \N \takeaway J_1}} \geq 2N.$
However, the completeness of $\Ac$
implies that the codimension of $\clspan\set{f_n}_{n \in \N \takeaway J_1}$
is at most $N,$ which is a contradiction.
\end{proof}

Now we prove that the excess of a frame equals the excess of
any alternative or synthesis pseudo-dual.

\begin{theorem} \label{excess_thm}
If $\Fc = \set{x_n}_{n \in \N}$ is a frame for $H$
then $e(\Fc) = e(\Yc) = e(\Zc)$
for every alternative dual $\Yc$ and synthesis pseudo-dual $\Zc.$
\end{theorem}
\begin{proof}
As we pointed out earlier, it suffices to consider the case
$e(\Fc) = \infty.$

\smallskip\emph{Alternative Duals}.
Let $\Yc = \set{y_n}_{n \in \N}$ be an alternative dual of $\Fc,$ and
suppose that $e(\Yc) = 0.$
In this case $\Yc$ is minimal, so there exists a biorthogonal sequence
$\set{\ty_n}_{n \in \N}.$
But since $\Yc$ is an alternative dual, this implies that
$\ty_k = \sum \, \ip{\ty_k}{y_n} \, x_n = x_k$ for every $k.$
Hence $\set{x_n}_{n \in \N}$ is minimal and so has zero excess,
which is a contradiction.

Next, suppose that $0 < e(\Yc) < \infty,$ and let $N = e(\Yc).$
By reindexing, we may assume $\set{y_n}_{n > N}$ has zero excess
and hence is minimal.
Consequently it has a biorthogonal sequence $\set{\ty_n}_{n > N}.$
Therefore, if $k > N$ then
$$\ty_k
\Eq \sum_{n=1}^\infty \, \ip{\ty_k}{y_n} \, x_n
\Eq x_k \Plus \sum_{n=1}^N \, \ip{\ty_k}{y_n} \, x_n.$$
For each $k > N$ let
$$w_k \Eq \sum_{n=1}^N \, \ip{\ty_k}{y_n} \, x_n,$$
so $x_k + w_k = \ty_k$ for $k > N.$
By Lemma \ref{minimal_lemma},
$$K
\Eq e\Bigparen{\set{w_n+x_n}_{n > N} \cup \set{x_n}_{n=1}^N}
\Lt \infty.$$
Since $\Fc$ has infinite excess, it is possible to remove $N+2K$ elements
from $\Fc$ yet leave the closed span unchanged.
Let $J_1 \subseteq \set{1,\dots,N}$ and $J_2 \subseteq \set{N+1,N+2,\dots}$
be such that $|J_1| + |J_2| = 2K+N$ and
$\clspan\set{x_n}_{n \notin J_1 \,\cup\, J_2} = H.$
Since $w_n \in \clspan\set{x_n}_{n=1}^N,$ it follows that
$$H
\Eq \clspan\bigparen{\set{x_n}_{n \in \set{1,\dots,N} \takeaway J_1} \cup
    \set{x_n}_{n > N,\, n \notin J_2}}
\Subseteq \clspan\bigparen{\set{x_n}_{n=1}^N \cup
    \set{x_n + w_n}_{n > N,\, n \notin J_2}}.$$
Thus it is possible to remove $|J_2|$ elements from
$\set{x_n}_{n=1}^N \cup \set{w_n + x_n}_{n > N}$
without changing its closed span.
But $|J_2| \geq 2K > K,$ so this is a contradiction.

We conclude that we must have $e(\Yc) = \infty.$

\smallskip\emph{Synthesis Pseudo-Duals}.
Let $\Zc = \set{z_n}_{n \in \N}$ be a synthesis pseudo-dual.
If $e(\Zc)=0,$ then it has a biorthogonal sequence $\set{\tz_n}_{n \in \N}.$
Consequently, if $x \in H$ then for every $k$ we have
$$\ip{x}{\tz_k}
\Eq \biggip{\sum_{n=1}^\infty \, \ip{x}{x_n} \,z_n}{\tz_k}
\Eq \ip{x}{x_k}.$$
This implies that $x_k = \tz_k$ for every $k,$ contradicting our
assumption that $e(\Fc) = \infty.$

Finally, suppose that $0 < e(\Zc) < \infty.$
We may assume then that $\set{z_n}_{n > N}$ is minimal and has
a biorthogonal sequence $\set{\tz_n}_{n > N}.$
For $k > N$ and $x \in H,$
\begin{align*}
\ip{x}{\tz_k}
\Eq \biggip{\sum_{n=1}^\infty \, \ip{x}{x_n} \, z_n}{\tz_k}
& \Eq \ip{x}{x_k} \Plus \biggip{\sum_{n=1}^N \, \ip{x}{x_n} \, z_n}{\tz_k}
      \\
& \Eq \biggip{x}{x_k \Plus \sum_{n=1}^N \,\ip{\tz_k}{z_n} \,x_n}.
\end{align*}
Consequently,
$\tz_k
= x_k \Plus \sum_{n=1}^N \, \ip{\tz_k}{z_n} \, x_n$ for $k > N.$
An argument similar to the one used for alternative duals
then leads to a contradiction.
\end{proof}

\section*{Acknowledgments}
We thank Hans Feichtinger for bringing reference \cite{SB20}
to our attention.

\end{document}